\documentclass{gtpart}

\title{Thick simplices and quasicategories}

\author{Ezra Getzler}
\givenname{Ezra}
\surname{Getzler}
\address{Department of Mathematics, Northwestern University, Evanston,
Illinois, USA}
\email{getzler@northwestern.edu}

\dedicatory{Dedicated to Andr\'e Joyal on the occasion of his 70th
  birthday}

\input diagxy

\def\hookrightarrow{\mathchoice
  {\DOTSB\lhook\joinrel\relbar\joinrel\rightarrow}
  {\DOTSB\lhook\joinrel\rightarrow}
  {\DOTSB\lhook\joinrel\rightarrow}
  {\DOTSB\lhook\joinrel\rightarrow}}

\def\to{\mathchoice
  {\DOTSB\relbar\joinrel\rightarrow}
  {\DOTSB\joinrel\rightarrow}
  {\DOTSB\joinrel\rightarrow}
  {\DOTSB\joinrel\rightarrow}}

\arxivreference{math.CT/1307.6924}
\arxivpassword{mx752}

\volumenumber{}
\issuenumber{}
\publicationyear{}
\papernumber{}
\startpage{}
\endpage{}
\doi{}
\MR{}
\Zbl{}
\received{}
\revised{}
\accepted{}
\published{}
\publishedonline{}
\proposed{}
\seconded{}
\corresponding{}
\editor{}
\version{}

\usepackage{amsrefs}

\usepackage[bbgreekl]{mathbbol}

\let\D=\undef
\let\L=\undef
\DeclareMathSymbol{\D}{\mathord}{AMSa}{"4D}
\DeclareMathSymbol{\L}{\mathord}{symbols}{"5E}
\DeclareMathSymbol{\LAMBDA}{\mathord}{bbold}{"03}
\DeclareMathSymbol{\DELTA}{\mathord}{bbold}{"01}
\numberwithin{equation}{section}

\newtheorem{theorem}{Theorem}
\newtheorem{proposition}[theorem]{Proposition}
\newtheorem{corollary}[theorem]{Corollary}

\theoremstyle{definition}

\newtheorem{definition}[theorem]{Definition}

\DeclareMathOperator{\cosk}{cosk}

\DeclareMathOperator{\Map}{Map}
\newcommand{\p}{\partial}
\newcommand{\eps}{\varepsilon}
\newcommand{\bull}{\bullet}
\newcommand{\CU}{\mathcal{U}}

\newcommand{\<}{\langle\,}
\renewcommand{\>}{\,\rangle}
\newcommand{\GG}{\mathbb{G}}
\newcommand{\G}{\mathsf{G}}
\newcommand{\deltabar}{\bar{\delta}}

\newcommand{\CG}{\mathcal{G}}

\begin{document}

\maketitle

\section{Thick simplices}

In the study of quasicategories, an important role is played by the
simplicial sets
\begin{equation*}
  \DELTA^n = \cosk_0\D^n .
\end{equation*}
In other words, $\DELTA^n$ is the simplicial set whose $k$-simplices
are the sequences
\begin{equation*}
  (i_0,\dots,i_k) \in \{0,\dots,n\}^{k+1} .
\end{equation*}
We call these simplicial sets the \textbf{thick simplices}. They go by
different names in the literature: Rezk~\cite{Rezk} denotes them
$E(n)$, while Joyal and Tierney~\cite{JT} use the notation $\D'[n]$.

The thick simplex $\DELTA^k$ is the nerve of the groupoid obtained by
localizing the category $[k]$, that is, adjoining inverses for all
non-identity morphisms. (This is the groupoid whose objects are the
vertices of $\D^k$, and such that the set of morphisms from $i$ to $j$
has exactly one element for each pair of vertices $i$ and $j$.)

A $k$-simplex of $\DELTA^n$ is non-degenerate if $i_{j-1}\ne i_j$ for
all $1\le j\le k$; $\DELTA^n$ has $(n+1)n^k$ non-degenerate
$k$-simplices. For $n=0$, $\DELTA^0=\D^0$.

The thick $1$-simplex $\DELTA^1$ has two non-degenerate simplices
$\delta_k=(0,1,0,1,\dots)$ and $\deltabar_k=(1,0,1,0,\dots)$ in each
dimension $k\ge0$, and $\DELTA^1$ may be identified with the
infinite-dimensional sphere $S^\infty$.

Recall that the horn $\L^n_i\subset\D^n$ is the simplicial set
obtained by removing the face $\p_i\D^n$ from the boundary $\p\D^n$ of
$\D^n$.

If $X_\bull$ is a simplicial set, form the coend
\begin{equation*}
  X_\bull\times_\D\DELTA^\bull = \textstyle\int^{n\in\D} X_n \times
  \DELTA^n .
\end{equation*}
(This simplicial set is denoted $k_!X_\bull$ by Joyal and Tierney
\cite{JT}.) In particular, we have the thick horns and the thick
boundary
\begin{equation*}
  \LAMBDA^n_i = \L^n_i\times_\D\DELTA^\bull \subset
  \p\DELTA^n = \p\D^n\times_\D\DELTA^\bull .
\end{equation*}
A simplex $(i_0,\dots,i_k)$ lies in $\p\DELTA^n$ if and only if
$\{i_0,\dots,i_k\}$ is a proper subset of the set $\{0,\dots,n\}$ of
vertices of $\DELTA^n$, and in $\LAMBDA^n_i$ if and only if
$\{i_0,\dots,i_k\}$ is a proper subset of the set $\{0,\dots,n\}$ of
vertices of $\DELTA^n$ not equal to
$\{0,\dots,\widehat{\imath},\dots,n\}$.

Consider the following cofibrations (inclusions) of simplicial sets:
\begin{enumerate}
\item the inclusion of the boundary
  \begin{equation*}
    \mu_n : \p\D^n \hookrightarrow \D^n ;
  \end{equation*}
\item a thickened form of this inclusion
  \begin{equation*}
    \bbmu_n : \p\DELTA^n \cup \D^n \hookrightarrow \D^n ;
  \end{equation*}
\item the inclusion of a horn
  \begin{equation*}
    \lambda_{n,i} : \L^n_i \hookrightarrow \D^n ;
  \end{equation*}
\item a thickened form of this inclusion
  \begin{equation*}
    \bblambda_{n,i} : \LAMBDA^n_i \cup \D^n \hookrightarrow \DELTA^n .
  \end{equation*}
\end{enumerate}

Recall the following result.

\begin{proposition}[Proposition~1.20, Joyal and Tierney~\cite{JT}]
  \label{JT}
  The inclusion
  \begin{equation*}
    \bbmu_n : \p\DELTA^n\cup\D^n \hookrightarrow \DELTA^n
  \end{equation*}
  is a trivial cofibration of simplicial sets.
\end{proposition}

The proof of this proposition comes down to the statement that the
geometric realization of $\p\DELTA^n\cup\D^n$ is contractible, and
ultimately, that $\DELTA^n$ is contractible for all $n\ge0$. This is
the case because $\DELTA^n$ is the nerve of a groupoid equivalent to
the trivial group.

Let $X_\bull$ be a simplicial set. The proposition immediately
implies, by induction on the simplices of $X_\bull$, that the
inclusion
\begin{equation*}
  X_\bull \hookrightarrow X_\bull\times_\D\DELTA^\bull
\end{equation*}
is a trivial cofibration.

If $\CU=\{U_\alpha\}$ is an open cover of a space $X$, and $A$ is a
sheaf of abelian groups on $X$, let $\check{C}^*(\CU,A)$ be the
\v{C}ech-complex
\begin{equation*}
  \check{C}^k(\CU,A) = \prod_{\alpha_0,\dots,\alpha_k}
  \Gamma(U_{\alpha_0}\cap\dots\cap U_{\alpha_k},A) ,
\end{equation*}
and let $\check{C}_\eps^*(\CU,A)\subset\check{C}^*(\CU,A)$ be the
subcomplex of antisymmetric \v{C}ech-cochains. Serre \cite{Serre}*{\S20,
  Proposition~2} proves that the inclusion of complexes
\begin{equation*}
  \check{C}_\eps^*(\CU,A) \hookrightarrow \check{C}^*(\CU,A)
\end{equation*}
is a quasi-isomorphism. Proposition~\ref{JT} may be viewed as a
non-abelian analogue of this result.

\section{Expansions and inner expansions}

Our goal is to refine Proposition~\ref{JT}, in such a way as to yield
more combinatorial proofs of certain basic results in the theory of
quasicategories.

\begin{definition}
  An \textbf{expansion} $S\hookrightarrow T$ is a map of simplicial
  sets such that there exists a filtration
  \begin{equation*}
    S=S_0 \subset S_1 \subset S_2 \subset \dots \subset T
  \end{equation*}
  satisfying the following conditions:
  \begin{enumerate}
  \item $T=\bigcup_kS_k$;
  \item there are maps $\L^{n_k}_{i_k} \to S_{k-1}$ such that the
    following diagram is a pushout square:
    \begin{equation*}
      \begin{xy}
        \Square[\L^{n_k}_{i_k}`S_{k-1}`\D^{n_k}`S_k;```]
      \end{xy}
    \end{equation*}
  \end{enumerate}
  The expansion is \textbf{inner} if $0<i_k<n_k$ for all $k$.
\end{definition}

\renewcommand{\Q}{\tilde{Q}}

Every expansion is a trivial cofibration. The first part of the
following proposition is thus a sharpening of Proposition~\ref{JT}.
\begin{proposition}
  \label{jt}
  \mbox{}
  \begin{enumerate}
  \item The inclusion
    \begin{equation*}
      \bbmu_n : \p\DELTA^n \cup \D^n \hookrightarrow \DELTA^n
    \end{equation*}
    is an expansion.
  \item If $0<i<n$, the inclusion
    \begin{equation*}
      \bblambda_{n,i} : \LAMBDA^n_i \cup \D^n \hookrightarrow \DELTA^n
    \end{equation*}
    is an inner expansion.
  \end{enumerate}
\end{proposition}
\begin{proof}
  We start with the proof of Part~a). Partition the non-degenerate
  simplices of $\DELTA^n$ into the following disjoint subsets: a
  $k$-simplex $s=(i_0\dots i_k)$ of $\DELTA^n$ has
  \begin{enumerate}
  \item type $P_k$ if it is contained in $\p\DELTA^n\cup\D^n$;
  \item type $Q_{k,m}$, $0\le m<n$, if there exists $m\ge0$ (which is
    necessarily uniquely determined by $s$) such that
    \begin{enumerate}
    \item $i_j=j$ for $0\le j\le m$, and
    \item $\{i_{m+1},\dots,i_n\}=\{m,\dots,n\}$;
    \end{enumerate}
  \item type $R_k$ otherwise.
  \end{enumerate}
  Write $Q_k$ for the unions of the types $Q_{k,m}$ over all values of
  $m$.

  We show that $\DELTA^n$ is obtained by adjoining the simplices of
  type $Q_{k,m}$ to $\p\DELTA^n\cup\D^n$ in order first of increasing
  $k$, then of decreasing $m$. (The order in which the simplices are
  adjoined within the sets $Q_{k,m}$ is unimportant.)

  Given a simplex $s=(i_0,\dots,i_k)$ of type $R_k$, let $m$ be
  the largest integer such that $i_j=j$ for $j<m$. Thus
  \begin{equation*}
    s = (0,\dots,m-1,i_m,\dots,i_k) ,
  \end{equation*}
  and $i_m\ne m$. The infimum $\ell$ of the set $\{i_m,\dots,i_k\}$
  equals $m$: it cannot be any larger, or the simplex would have type
  $P_k$, and it cannot be any smaller, or the simplex would have type
  $Q_{k,\ell}$. Define the simplex
  \begin{equation*}
    s' = (0,\dots,m,i_m,\dots,i_k)
  \end{equation*}
  of type $Q_{k+1,m}$. We have $s=\p_ms'$.

  If $m$ occurs more than once in the sequence $\{i_m,\dots,i_k\}$,
  then all faces of the simplex $s'$ are either degeneracices of
  $(k-1)$-simplices, or of type $P_k$ or $Q_k$. If $m$ does occur just
  once in this sequence, say $i_\ell=m$, then it is still the case
  that all faces of of the simplex $s'$ other than $\p_ms'$ and
  $\p_{\ell+1}s'$ are either degeneracices of $(k-1)$-simplices, or of
  type $P_k$ or $Q_k$. The face
  \begin{equation*}
    \p_{\ell+1}s' =
    (0,\dots,m,i_m,\dots,\widehat{\imath}_\ell,\dots,i_k)
  \end{equation*}
  is of type $R_k$: it is a face of a simplex of type $Q_{k+1,m'}$,
  where $m'>m$.

  In any case, all faces of $s'$ are simplices of type $P_j$ or $Q_j$,
  $j\le k$, and simplices of type $Q_{k+1,m'}$, where $m'>m$. This
  completes the proof of Part~a).

  We now turn to the proof of Part b). Partition the non-degenerate
  simplices of $\DELTA^n$ into the following disjoint subsets: a
  $k$-simplex $s=(i_0\dots i_k)$ of $\DELTA^n$ has
  \begin{enumerate}
  \item type $P_k$ if it is contained in $\LAMBDA^n_i\cup\D^n$;
  \item type $Q_{k,m}$, $m\ge0$ even, if it is not of type $P_k$, has
    the form
    \begin{equation*}
      (i_0,\dots,i_{j-2},\underbrace{j,j',j,j',\dots,j}_{\text{length $m+3$}},
      i_{j+m+2},\dots,i_k) ,
    \end{equation*}
    and $i_{j+m+2}$ does not equal $j'$;
  \item type $\Q_{k,m}$, $m\ge0$, if it is not of type $P_k$, has the
    form
    \begin{equation*}
      (i_0,\dots,i_{j-1},\underbrace{j,j',j,j',\dots}_{\text{length $m+2$}},
      i_{j+m+2},\dots,i_k) ,
    \end{equation*}
    and $i_{j+m+1}$ equals neither $j$ nor $j'$;
  \item type $R_k$ otherwise.
  \end{enumerate}
  Write $Q_k$ and $\Q_k$ for the unions of the types $Q_{k,m}$ and
  $\Q_{k,m}$ over all values of $m$. We show that $\DELTA^n$ is
  obtained by adjoining the simplices of type $Q_k$ and $\Q_k$ in
  order of increasing $k$, with the simplices of $Q_k$ adjoined before
  those of $\Q_k$. (The order in which the simplices are adjoined
  within the sets $Q_k$ and $\Q_k$ is unimportant.)

  Let $s=(i_0,\dots,i_k)$ be a simplex of type $R_k$, and let
  $j'=i_j$. We will assign a simplex $s'$ of type $Q_{k+1}$ or
  $\Q_{k+1}$ such that $s$ is an inner face of $s'$.

  First consider the case in which $i_{j-1}=j$. Let $m$ be the largest
  even number such that $s$ has the form
  \begin{equation*}
    (i_0,\dots,i_{j-2},\underbrace{j,j',j,j',\dots}_{\text{length $m+2$}},
    i_{j+m+1},\dots,i_k) ,
  \end{equation*}
  and define $s'$ to be the simplex
  \begin{equation*}
    s' = (i_0,\dots,i_{j-2},\underbrace{j,j',j,j',\dots,j}_{\text{length $m+3$}},
    i_{j+m+1},\dots,i_k)
  \end{equation*}
  of type $Q_{k+1,m}$. The condition $0<j<n$ implies that
  $s=\p_{j+m+1}s'$ is an inner face of $s'$: otherwise we would
  have $k=j+m$, and $\{i_0,\dots,i_k\}=\{i_0,\dots,i_j\}$, which is
  a proper subset of $\{0,\dots,n\}$.

  If $i_{j-1}\ne j$, we let $m$ be the largest natural number such
  that $s$ has the form
  \begin{equation*}
    (i_0,\dots,i_{j-1},\underbrace{j',j,j',j,\dots}_{\text{length $m+1$}},
    i_{j+m+1},\dots,i_k) ,
  \end{equation*}
  and define $s'$ to be the simplex
  \begin{equation*}
    s' = (i_0,\dots,i_{j-1},\underbrace{j,j',j,j',\dots}_{\text{length $m+2$}},
    i_{j+m+1},\dots,i_k)
  \end{equation*}
  of type $\Q_{k+1,m}$. The condition $0<j<n$ again implies that
  $s=\p_js'$ is an inner face of $s'$.

  Every face of a simplex $s'$ of type $Q_{k+1}$ other than
  $\p_{k+m+1}s'$, is either a degeneracy of a $(k-1)$-simplex, or of
  type $P_k$, type $Q_k$ or type $\Q_k$.

  Likewise, every face of a simplex $s'$ of type $\Q_{k+1,m}$, $m>0$,
  other than $\p_js'$, is either a degeneracy of a $(k-1)$-simplex, or
  of type $P_k$, type $Q_k$ or type $\Q_k$.

  It remains to consider a simplex $s'$ of type $\Q_{k+1,0}$. The face
  $\p_\ell s'$, $\ell>j$, is either a degeneracy of a $(k-1)$-simplex,
  or of type $P_k$, type $Q_k$ or type $\Q_k$. The face $\p_\ell s'$,
  $\ell<j$, is either a degeneracy of a $(k-1)$-simplex, of type
  $P_k$, or the face $\p_{j+1}s''$ of the simplex
  \begin{equation*}
    s'' = (i_0,\dots,\widehat{\imath}_\ell,\dots,i_{j-1},j,j',j,i_{j+2},\dots,i_k)
  \end{equation*}
  of type $Q_{k+1,0}$. This completes the proof of the proposition.
\end{proof}

\begin{corollary}
  \label{lambda}
  If $n>1$, the inclusion $\LAMBDA^n_i \hookrightarrow \DELTA^n$ is an
  inner expansion for all $0\le i\le n$.
\end{corollary}
\begin{proof}
  The action of the symmetric group $S_{n+1}$ on the simplicial set
  $\DELTA^n$ induces a transitive permutation of the subcomplexes
  $\L^n_i$. Thus, it suffices to establish the result when $0<i<n$, in
  which case we may factor the inclusion $\LAMBDA^n_i \hookrightarrow
  \DELTA^n$ as the inner expansion
  \begin{equation*}
    \LAMBDA^n_i \hookrightarrow \LAMBDA^n_i \cup \D^n
  \end{equation*}
  followed by the inner expansion $\bblambda_{n,i} : \LAMBDA^n_i \cup
  \D^n \hookrightarrow \DELTA^n$.
\end{proof}

\section{A K\"unneth theorem}

Let $\D^{m,n}=\D^m\times\D^n$ be the \textbf{prism}. There are two
different thick prisms. The first of these is the product of two thick
simplices, which is itself a thick simplex:
\begin{equation*}
  \DELTA^m \times \DELTA^n \cong \DELTA^{mn+m+n} .
\end{equation*}
In fact, the vertices of the product are the pairs $(i,j)$, where
$0\le i\le m$ is a vertex of $\DELTA^m$, and $0\le j\le n$ is a vertex
of $\DELTA^n$. This simplicial set contains the simplicial set
\begin{equation*}
  \DELTA^{m,n} = \D^{m,n} \times_\D \DELTA^\bull ,
\end{equation*}
which is the union of $\binom{m+n}{m}$ copies of the thick simplex
$\DELTA^{m+n}$. The following proposition is an analogue of the
K\"unneth theorem.
\begin{proposition}
  \label{Kunneth}
  The inclusion $\DELTA^{m,n}\hookrightarrow\DELTA^m\times\DELTA^n$ is
  an inner expansion.
\end{proposition}
\begin{proof}
  If $S\hookrightarrow T$ is an expansion of simplicial sets such that
  all of the simplices adjoined to $S$ in order to obtain $T$ have
  dimension greater than $1$, then by Corollary~\ref{lambda},
  \begin{equation*}
    S\times_\D\DELTA^\bull \hookrightarrow T\times_\D\DELTA^\bull
  \end{equation*}
  is an inner expansion. Thus, to prove the proposition, it suffices
  to prove that the inclusion
  \begin{equation*}
    \D^{m,n} \hookrightarrow \D^{mn+m+n}
  \end{equation*}
  is an expansion.

  Partition the non-degenerate $k$-simplices of $\D^{mn+m+n}$ into the
  following disjoint subsets: a $k$-simplex
  \begin{align*}
    s &= ((i_0,j_0),\dots,(i_k,j_k)) \\
    &= (i_0,\dots,i_k) \times (j_0,\dots,j_k)
  \end{align*}
  of $\D^{mn+m+n}$ has
  \begin{enumerate}
  \item type $P_k$ if it is contained in $\D^{m,n}$;
  \item type $Q_k$ if it is not of type $P_k$, and $i_0=j_0=0$;
  \item type $R_k$ otherwise.
  \end{enumerate}
  We show that $\D^{mn+m+n}$ is obtained from $\D^{m,n}$ by adjoining
  the simplices of type $Q_k$ in order of increasing $k$. (The order
  in which the simplices are adjoined within the sets $Q_k$ is
  unimportant.)

  If $s=((i_0,j_0),\dots,(i_k,j_k))$ is a simplex of type $R_k$, let
  $s'$ be the simplex
  \begin{align*}
    s' &= ((0,0),(i_0,j_0),\dots,(i_k,j_k)) \\
    &= (0,i_0,\dots,i_k) \times (0,j_0,\dots,j_k)
  \end{align*}
  of type $Q_{k+1}$. Thus $s=\p_0s'$, while every other face of $s'$ has
  type either $P_k$ or $Q_k$. This completes the proof.
\end{proof}

\section{Application to the theory of quasicategories}

If $S_\bull$ and $X_\bull$ are simplicial sets, let $\Map(S,X)$ be the
set of simplicial maps from $S_\bull$ to $X_\bull$. If
$i:S\hookrightarrow T$ is a cofibration of simplicial sets (that is,
if $S_k\to T_k$ is a monomorphism for all $k\ge0$), let $\<i\mid f\>$
be the induced function from $\Map(T,X)$ to
\begin{equation*}
  \Map(S,X) \times_{\Map(T,X)} \Map(T,Y) .
\end{equation*}
We write $\<i\mid X\>$ in place of $\<i\mid f\>$, when $f$ is the
unique simplicial map from the simplicial set $X_\bull$ to the
terminal simplicial set $\Delta^0$.

\begin{definition}
  A simplicial set $X_\bull$ is
  \begin{enumerate}
  \item a \textbf{quasicategory} if the functions
    $\<\lambda_{n,i}\mid X\>$ are surjective for $0<i<n$;
  \item a \textbf{Kan complex} if the functions $\<\lambda_{n,i}\mid
    X\>$ are surjective for $n>0$ and $0\le i\le n$.
  \end{enumerate}
  A simplicial map $f:X_\bull\to Y_\bull$ is a \textbf{trivial
    fibration} if the functions $\<\mu_n\mid f\>$ are surjective for
  $n\ge0$.
\end{definition}

If $X_\bull$ is a simplicial set, denote by $\GG(X)_\bull$ the
simplicial set with $n$-simplices
\begin{equation*}
  \GG(X)_n = \Map(\DELTA^n,X) .
\end{equation*}
The inclusions $\D^n\hookrightarrow\DELTA^n$ induce a natural
simplicial map from $\GG(X)_\bull$ to $X_\bull$, with image
$\G(X)_\bull$. A functor similar to $\GG$ was introduced by
Rezk~\cite{Rezk}, in his study of complete Segal spaces.

\begin{proposition}
  \label{GGG}
  $\GG(\GG(X))_\bull\cong\GG(\G(X))_\bull\cong\GG(X)_\bull$
\end{proposition}
\begin{proof}
  In order to prove that $\GG(\GG(X))_\bull$ is isomorphic to
  $\GG(X)_\bull$, it suffices to show that for all $k,n\ge0$,
  \begin{equation*}
    \Map(\D^k,\DELTA^n)\cong\Map(\DELTA^k,\DELTA^n) .
  \end{equation*}
  Since $\DELTA^k$ is the nerve of the groupoid obtained by localizing
  the category $[k]$, we see that $\Map(\DELTA^k,\DELTA^n)$ may be
  identified with the set of functors from the localization of $[k]$
  to the localization of $[n]$. But a functor from the localization of
  $[k]$ to a groupoid $\CG$ determines, and is determined by, a
  functor from $[k]$ to $\CG$, i.e.\ by a $k$-simplex of the nerve
  $N_\bull\CG$ of $\CG$.

  Applying the functor $\GG_n$ to the composition of morphisms
  \begin{equation*}
    \GG(X)_\bull \to \G(X)_\bull \to X_\bull ,
  \end{equation*}
  we obtain a factorization of the identity map of $\GG(X)_n$:
  \begin{equation*}
    \GG(\GG(X))_n \cong \GG(X)_n \to \GG(\G(X))_n \to \GG(X)_n .
  \end{equation*}
  Since the functor $\GG_n$ is a finite limit, it takes monomorphisms
  to monomorphisms. Thus the morphism from $\GG(\G(X))_n$ to
  $\GG(X)_n$ is a monomorphism, and hence, since it has a section, an
  isomorphism.
\end{proof}

The following theorem assembles results of Joyal \cite{Joyal}, Joyal
and Tierney~\cite{JT} and Rezk~\cite{Rezk}. Our approach is to reduce
the proof to Proposition~\ref{jt}. Although these results are not new
in the context of quasicategories, the approach developed here allows
the development of a theory of quasicategories in the setting where
sets are replaced by analytic spaces (or \textbf{Lie
  quasicategories}, by analogy with Lie groupoids, of which they are
a generalization).
\begin{theorem}
  \label{main}
  Let $X_\bull$ be a quasicategory.
  \begin{enumerate}
  \item The simplicial sets $\GG(X)_\bull$ and $\G(X)_\bull$ are Kan
    complexes.
  \item The simplicial map $\GG(X)_\bull\to \G(X)_\bull$ is a trivial
    fibration.
  \end{enumerate}
\end{theorem}
\begin{proof}
  The simplicial set $\GG(X)_\bull$ is a Kan complex by
  Corollary~\ref{lambda}.

  For $n>0$, consider the assertions \\
  \noindent {\boldmath$\textup{A}_n$:} for all $0\le i\le n$, the function
  $\G(X)_n\to\L_{n,i}(\G(X))$ is surjective; and \\
  \noindent {\boldmath$\textup{B}_n$:} for all $0\le i\le n$, the function
  \begin{equation*}
    \GG(X)_n \to \L_{n,i}(\GG(X))\times_{\L_{n,i}(\G(X))}\G(X)_n
  \end{equation*}
  is surjective.

  We will establish these assertions for all $n$ by induction: it will
  follow that $\G(X)_\bull$ is a Kan complex.

  $\textup{A}_1$ is clear. If $\textup{A}_n$ holds, the composition
  \begin{equation*}
    \GG(X)_n\times_{\L_{n,i}(\G(X))}\G(X)_n \to
    \L_{n,i}(\GG(X))\times_{\L_{n,i}(\G(X))}\G(X)_n \to \GG(X)_n
  \end{equation*}
  is surjective, establishing $\textup{B}_n$.

  Suppose that $T$ is a finite simplicial set and $S\hookrightarrow T$
  is an expansion obtained by attaching simplices of dimension at most
  $n-1$ to $S$. If $\textup{B}_{n-1}$ holds, the function
  $\<S\hookrightarrow T\mid\GG(X)\to\G(X)\>$ is surjective. In
  particular, applying this argument to the expansion
  $\D^0\hookrightarrow\L^n_i$, we see that the function
  \begin{equation*}
    \L_{n,i}(\GG(X)) \to \GG(X)_0 \times_{\G(X)_0} \L_{n,i}(\G(X))
    \cong \L_{n,i}(\G(X))
  \end{equation*}
  is surjective. In the commuting diagram
  \begin{equation*}
    \begin{xy}
      \Square/{>}`{.>}`{>}`{.>}/%
      [\GG(X)_n`\L_{n,i}(\GG(X))`\G(X)_n`\L_{n,i}(\G(X));```]      
    \end{xy}
  \end{equation*}
  the solid arrows are surjective, hence so is the bottom arrow,
  establishing $\textup{A}_n$, and completing the induction.

  We may now prove that $\GG(X)_\bull\to\G(X)_\bull$ is a trivial
  fibration. We must show that for each $n\ge0$, the function
  \begin{equation*}
    \GG(X)_n \to M_n(\GG(X)_n) \times_{M_n(\G(X))_n} \G(X)_n
  \end{equation*}
  is a surjection. By Proposition~\ref{GGG}, this function may be
  identified with the function
  \begin{equation*}
    \GG(\G(X))_n \to M_n(\GG(\G(X))_n) \times_{M_n(\G(X))_n} \G(X)_n ,
  \end{equation*}
  and hence with $\<\bbmu_n\mid\G(X)\>$. As $\bbmu_n$ is an expansion
  and $\G(X)_\bull$ is a Kan complex, the result follows.
\end{proof}

\paragraph*{Acknowledgements}

The results of this paper were motivated by joint work with Kai
Behrend on the theory of higher stacks.

This paper was conceived while the author held a Microsoft Research
Visiting Fellowship at the Isaac Newton Institute of Mathematical
Sciences of Cambridge University. He thanks the organizers of the
programme on Operads and Multiple Zeta-Values, John D.S. Jones and
Bruno Vallette, for this opportunity, the other participants in the
programme, and particularly Mikhail Kapranov, Ieke Moerdijk and Jon
Pridham, for useful conversations on the subject of this paper, and
the Simons Foundation for support under a Collaboration Grant for
Mathematicians.


\begin{bibdiv}

\begin{biblist}

\bib{Joyal}{article}{
   author={Joyal, A.},
   title={Quasi-categories and Kan complexes},
   journal={J. Pure Appl. Algebra},
   volume={175},
   date={2002},
   number={1-3},
   pages={207--222},
}


\bib{JT}{article}{
   author={Joyal, Andr{\'e}},
   author={Tierney, Myles},
   title={Quasi-categories vs Segal spaces},
   conference={
      title={Categories in algebra, geometry and mathematical physics},
   },
   book={
      series={Contemp. Math.},
      volume={431},
      publisher={Amer. Math. Soc.},
      place={Providence, RI},
   },
   date={2007},
   pages={277--326},
}

\bib{Rezk}{article}{
   author={Rezk, Charles},
   title={A model for the homotopy theory of homotopy theory},
   journal={Trans. Amer. Math. Soc.},
   volume={353},
   date={2001},
   number={3},
   pages={973--1007},
}

\bib{Serre}{article}{
   author={Serre, Jean-Pierre},
   title={Faisceaux alg\'ebriques coh\'erents},
   journal={Ann. of Math. (2)},
   volume={61},
   date={1955},
   pages={197--278},
}

\end{biblist}

\end{bibdiv}
\end{document}